\documentclass[draft]{amsart}

\usepackage{amsmath}
\usepackage{amssymb}
\usepackage{url}

\newtheorem{thm}{Theorem}[section]
\newtheorem{cor}[thm]{Corollary}
\newtheorem{lem}[thm]{Lemma}
\newtheorem{prop}[thm]{Proposition}

\theoremstyle{definition}
\newtheorem*{defn}{Definition}
\newtheorem{assumption}{Assumption} 
\newtheorem*{example}{Example}

\theoremstyle{remark}
\newtheorem*{rem}{Remark}

\numberwithin{equation}{section}

\newcommand{\C}{C}
\newcommand{\fq}{\mathbb{F}_{q}}
\newcommand{\fp}{\mathbb{F}_{p}}
\newcommand{\fqk}{\mathbb{F}_{q^k}}
\newcommand\Z{{\mathbb Z}}

\newcommand\F{{\mathbb F}}
\newcommand\Q{{\mathbb Q}}
\renewcommand\O{{\mathcal O}}
\newcommand\Gal{\mbox{\rm Gal}}
\newcommand\End{\mbox{\rm End}}

\newcommand\Jac{\mbox{\rm Jac}}
\newcommand\Tr{\mbox{\rm Tr}}

\begin{document}

\title{Distortion maps for genus two curves}

\author[Galbraith]{Steven D.~Galbraith}
\address{Mathematics Department,
    Royal Holloway University of London,
    Egham, Surrey TW20 0EX, United Kingdom.
}
\email{Steven.Galbraith@rhul.ac.uk}
\urladdr{http://www.isg.rhul.ac.uk/~sdg/}

\author[Pujol\`as]{Jordi Pujol\`as}
\address{Departament de Matem\`atica Aplicada 4,
    Universitat Polit\`ecnica de Catalunya,
    Jordi Girona 1-3, 08034 Barcelona, Spain.
}\email{jpujolas@matematica.udl.es}

\author[Ritzenthaler]{Christophe Ritzenthaler}
\address{Institut de Math\'ematiques de Luminy,
         UMR 6206 du CNRS,
         Luminy, Case 907, 13288 Marseille, France.}
\email{ritzenth@iml.univ-mrs.fr}
\urladdr{http://iml.univ-mrs.fr/~ritzenth/}

\author[Smith]{Benjamin Smith}
\address{Mathematics Department,
    Royal Holloway University of London,
    Egham, Surrey TW20 0EX, United Kingdom.
}
\email{Ben.Smith@rhul.ac.uk}
\urladdr{http://www.ma.rhul.ac.uk/~bensmith/}

\begin{abstract}
Distortion maps are a useful tool for pairing based cryptography.
Compared with elliptic curves, the case of hyperelliptic curves
of genus $g > 1$ is more complicated since the full torsion subgroup
has rank $2g$.
In this paper we prove that
distortion maps always exist for supersingular curves of
genus $g>1$ and we give several examples in genus $2$.

\end{abstract}

\maketitle


\section{Introduction}
\label{sec:intro}

Let $q$ be a power of a prime $p$.
Let $\C$ be a
non-singular, geometrically irreducible, and projective
curve
defined over the finite field $\fq$.
The Jacobian variety of $C$ is denoted by $\Jac(\C)$,
and the $q$-power Frobenius map is denoted $\pi$.
Throughout this paper, we identify $\Jac (\C)(\F_{q^n})$ with the
degree zero divisor class group of $\C$ over $\F_{q^n}$.
Let $r$ be a prime number dividing $\# \Jac (\C)(\fq)$
and coprime to $p$.
We define the \emph{embedding degree} to be the smallest positive integer $k$
such that $r$ divides $q^k-1$.
Note that $\fqk$ is then
the field generated over $\fq$
by adjoining the $r^\mathrm{th}$ roots of unity.
If $A$ is an abelian variety,
then $\End_K(A)$
denotes the ring of endomorphisms of $A$ defined over a field $K$,
and $\End(A)$ the ring of endomorphisms of $A$
defined over an algebraic closure of $K$.
Unless specified otherwise,
all morphisms are defined over the algebraic closure of the field.

An elliptic curve $E$ over $\fq$ is called \emph{supersingular}
if the number of points on $E$ over $\fq$ is congruent to $1$ modulo $p$.
If $E$ is a supersingular elliptic curve,
then $\End( E )$ is an order in a quaternion algebra.
More generally,
an abelian variety $A$ of dimension $g$ over $\fq$ is called supersingular
if $A$ is isogenous over $\overline{\F}_q$ to a product $E^g$,
where $E$ is a supersingular elliptic curve.
In this case,
it follows that
$\End^0(A) = \End( A ) \otimes_{\Z} \Q$ is a $\Q$-algebra
of dimension $(2g)^2$ as a $\Q$-vector space.
Finally, a curve $\C$ is called supersingular
if $\Jac( \C )$ is a supersingular abelian variety.

Let $r>2$ be a prime dividing $q^k - 1$, and coprime to $q$. The
Tate pairing (see Frey and R\"uck \cite{FreyRueck}) is a
non-degenerate bilinear pairing of the $r$-torsion in the divisor
class group of $\C$ over $\fqk$ with a certain quotient group of the
divisor class group over $\fqk$. Using standard methods (see
\cite{BSS2,Gal}), we can obtain from the Tate pairing a bilinear
pairing (often called the reduced Tate pairing) $e_r( \cdot, \cdot)$
from $\Jac( \C )(\fqk)[r]$ to the group $\mu_r$ of $r^\mathrm{th}$
roots of unity in $\fqk^*$.

When $\C$ is supersingular
and $r \Vert \#\Jac( \C )( \fq )$,
then $\Jac(\C)[r]$ is often contained in $\Jac( \C )( \fqk )$
(see \cite{S-X}).
In this case,
the Weil pairing is also a non-degenerate bilinear
pairing on $\Jac( \C )(\fqk)[r]$.
If the embedding degree $k$ is small,
then either the Weil or Tate pairing may be useful for implementing
pairing-based cryptosystems
(see \cite{Gagne,Kenny1,Kenny2} for a survey).
We use the notation $e_r(\cdot,\cdot)$
to denote any bilinear, nondegenerate, Galois-invariant pairing
on $\Jac(C)[r]$
(such as the Weil or reduced Tate pairings).

When $\Jac (\C)$ is supersingular,
the embedding degree $k$ is known to be
bounded above by a constant $k(g)$ depending only on the genus $g$ of $C$
(see \cite{Gal,RS}).
In cryptographic applications,
one tends to be interested in cases where
the embedding degree is greater than $1$,
but not ``too large''.

Bilinearity is an important property of pairings in cryptography:
for all integers $a$ and $b$ and elements $D_1$ and $D_2$ of $\Jac(C)[r]$,
we have $e_r( aD_1, bD_2 ) = e_r( D_1, D_2)^{ab}$.
For bilinearity to be useful, however,
it is necessary that
$e_r( D_1, D_2 ) \ne 1$.
It is known
that the Weil and Tate pairings are non-degenerate:
that is,
for each non-zero divisor class $D_1$ of order $r$,
there is a divisor class $D_2$ such that $e_r( D_1, D_2 ) \ne 1$.
A problem arises
when one wants to pair two specific divisors $D_1$ and $D_2$
such that $e_r( D_1, D_2 ) = 1$
--- this can happen, for example,
when
for efficiency reasons both divisors are defined over $\fq$, and $k>1$.
In these cases, we need distortion maps.

\begin{defn}
    A \emph{distortion map}
    for a non-degenerate pairing $e_r$
    and non-zero divisor classes $D_1$, $D_2$ of prime order  $r$
    on $\C$ is an endomorphism $\psi$ of $\Jac
    (\C)$ such that $e_r(D_1, \psi(D_2)) \neq 1$.
\end{defn}

Distortion maps were introduced by Verheul \cite{Verheul1}
for elliptic curves in the case where
$D_1$ and $D_2$ are defined over the ground field.
We stress that our definition depends on the choice of
divisor classes (and also the pairing).
In general,
it is not true
that there is a single choice of $\psi$
that is a distortion map for all pairs 
of non-zero divisor classes.

The goal of this paper
is to provide, for certain curves, a collection of
efficiently computable endomorphisms
such that
there is a suitable distortion map in the collection
for any pair 
of divisor classes on the curves.
Note that the Frobenius or trace maps may be used
as distortion maps in many situations,
including the case of ordinary curves;
but distortion maps for every pair 
can {\em only} be obtained for supersingular curves.

The case where $C$ is an elliptic curve is quite simple.
If $D_1$ and $D_2$ are nonzero divisor classes
and
$e_r( D_1, D_2 ) = 1$,
then any divisor $D_3$ of order $r$
which is independent of $D_2$
(that is, $\langle D_2 \rangle \cap \langle D_3 \rangle = \{ 0 \}$)
satisfies $e_r( D_1, D_3 ) \neq 1$.
This follows from the
non-degeneracy of the pairing,
and the fact that
the $r$-torsion of an elliptic curve has rank $2$.
For this reason, and others,
the problem of finding distortion maps for elliptic curves
is relatively easy to handle.
An algorithm to find distortion maps
for any supersingular elliptic curve has been given by Galbraith
and Rotger~\cite{GR}.\footnote{
    Note that there is a missing condition in Lemma~5.1 of \cite{GR},
    namely that $\psi(P) \ne 0$.
    Since the degree of $\psi$ in \cite{GR} is $d$,
    which is much smaller than $r$,
    this condition is always satisfied.
}

For curves $C$ of genus $g>1$,
the $r$-torsion of the Jacobian has rank $2g$;
so independence of divisors is not sufficient to imply
non-triviality of their pairing.
Indeed,
elementary linear algebra implies that
for every non-trivial divisor $D$ of order $r$,
there exists a basis for $\Jac(C)[r]$
such that $D$ pairs trivially with all but one of the basis elements.
Furthermore,
elements of $\End( \Jac( \C ))$ may be difficult to handle,
as they generally do not correspond to maps from $\C$ to itself.

In this paper, we discuss this situation, with particular emphasis
on curves of genus $2$. In Section~\ref{exist}, we prove that
distortion maps always exist for supersingular abelian varieties.
The rest of the paper is concerned with the question of whether such
maps can be easily computed on Jacobians of curves. In
Section~\ref{sec:curvelist}, we provide a list of examples of
supersingular curves with suitable embedding degree. These examples
are explored in depth in the subsequent sections. The results of
Section~\ref{Sec-3}, first presented in \cite{GP}, concern the case
$k=4$ where $p \equiv 2,3 \pmod{5}$. After illustrating our approach
on this simple case, we generalize the method to other curves in
Section~\ref{Sec-VW}. Section~\ref{sec:k5} deals with the case $k=5$
when $p=5$, and Section~\ref{sec:twist} deals with the case $k=6$
when $p \equiv 2 \pmod{3}$ (and $p \ne 2$). Finally,
Section~\ref{sec:char2} treats the case $k=12$ and $p=2$. We provide
non-trivial, efficient, explicit distortion maps for each curve.

Note that there is an important distinction between the cases
$k=4,12$ and $k=5,6$.
For the former cases,
our result are conditional,
since they depend on the assumption (verified in practice) that some
denominators can be canceled, which is the case if they are prime to $r$.
In these cases, there seems to be no easy explicit decomposition of
$\Jac(C)$: even when we know that the Jacobian splits into a product
$E \times E$, the degree of the induced morphisms from $C$ to $E$ is
unknown. However, the curves we consider in the cases $k=5$ and
$k=6$ are both twists of $y^2=x^6+1$, which has two degree-$2$ maps
to an elliptic curve $E$ (see Section~\ref{sec:twist} for details).
This structure is used in a crucial way to remove the assumption on
the denominators for the case~$k=6$ (for the case $k=5$ another
argument is used, which is restricted to the case where $e_r$ is the
Tate pairing
--- but the same proof could be adapted).

\section{The existence of distortion maps}
\label{exist}

Schoof and Verheul \cite{Verheul2} have shown
that distortion maps always exist for supersingular elliptic curves over $\fq$.
In this section,
we generalise their result to supersingular abelian varieties.

First, we recall an important theorem of Tate \cite{Tate}.
Suppose $A$ is an abelian
variety over a finite field $K$ of characteristic $p$,
and let $G  = \Gal( \overline{K} / K )$.
Let $l$ be a prime not equal to $p$,
and let $T_l(A) := \varprojlim A[l^n]$ be the $l$-Tate module of $A$.
Let $\End_G ( T_l( A ))$ denote
the ring of endomorphisms of $T_l(A)$
which commute with the action of $G$.
Tate's theorem states that the canonical injection
\[
   \End_{K} (A) \otimes_\Z \Z_l \longrightarrow \End_G( T_l(A))
\]
is an isomorphism.

\begin{thm}
\label{thm-exist}
    Let $A$ be a supersingular abelian variety of dimension $g$ over $\fq$,
    and let $r$ be a prime not equal to the characteristic of $\fq$.
    For every two non-trivial elements $D_1$ and $D_2$ of $A(\fq)[r]$,
    there exists an endomorphism $\phi$ of $A$
    such that $e_r( D_1, \phi( D_2 )) \neq 1$.
\end{thm}

\begin{proof}
    Let $d$ be an integer
    such that the $q^d$-power Frobenius map acts as
    an integer multiplication on $A$.
    Let $K = \F_{q^d}$ and $G = \Gal( \overline{K}/K )$.
    Since $A$ is supersingular,
    $\End( A ) \otimes_{\Z} \Z_r$
    is a free $\Z_r$-module of rank $(2g)^2$.
    By definition,
    $\End_K( A )$ is contained in $\End( A )$,
    so we may view $\End_K( A ) \otimes_{\Z} \Z_r$
    as a submodule of $\End( A ) \otimes_{\Z} \Z_r$.
    By Tate's theorem,
    $\End_K( A ) \otimes_{\Z} \Z_r$
    is isomorphic to the $\Z_r$-module
    $\End_G( T_r( A ))$
    of endomorphisms which commute with the $q^d$-power Frobenius
    --- but the $q^d$-power Frobenius is an integer,
    so it commutes with every endomorphism of $A$ (and $T_r(A)$).
    Thus $\End_G(T_r(A)) = \End(T_r(A))$.
    Since $T_r( A ) \cong \Z_r^{2g}$ as
    a $\Z_r$-module, we have
    $$
        \End_K(A)\otimes_{\Z}\Z_r
        \cong
        \End_G ( T_r( A ))
        \cong
        M_{2g}( \Z_r ).
    $$
    Hence $\End_K( A ) \cong \End( A )$
    also has rank $(2g)^2$.
    By restriction, we have
    $$
        \End_K ( A ) \otimes_{\Z} \Z/r\Z \cong M_{2g}( \Z/r\Z ).
    $$

    Let $D_3$ be an element of $A[r]$ such that $e_r( D_1, D_3 ) \neq 1$
    (in fact, $D_3$ is $K$-rational).
    There exists some matrix $\Phi$ in $M_{2g}( \Z / r \Z )$
    corresponding to a mapping of the subspace
    $\langle D_2 \rangle$ to $\langle D_3 \rangle$.
    Let $\phi$ be a preimage in $\End( A )$ of $\Phi$:
    by construction, $e_r( D_1, \phi( D_2 )) \neq 1$.
\end{proof}

The proof of Theorem~\ref{thm-exist}
shows that to have a distortion map for every pair of divisors,
we must have a full rank-$(2g)^2$ module of endomorphisms.
In other words,
if the rank of $\End(\Jac( \C ))$
is strictly less than $(2g)^2$,
then there will exist
non-zero elements
$D_1$ and $D_2$ of $\Jac( \C )[r]$ such that
$e_r( D_1, \psi( D_2 )) = 1$
for every endomorphism $\psi$ of $\Jac( \C )$.
In particular,
if $\C$ is not supersingular,
then there are pairs $(D_1,D_2)$
for which no distortion maps exist.

\begin{rem}
    It is important to note that Theorem~\ref{thm-exist}
    is not constructive.
\end{rem}

\section{Embedding degrees of supersingular genus~$2$ curves}
\label{sec:curvelist}

In this section we list some supersingular genus~$2$ curves which
are of potential interest for applications.
First, we recall the results of Rubin and Silverberg \cite{RS}
classifying the possible embedding
degrees for supersingular abelian varieties of dimension~$2$.
We focus on the case where $q$ is an odd power of $p$:
this gives the largest values for $k$,
and so is usually the most interesting case in practice.

\begin{thm}[Rubin--Silverberg~\cite{RS}]
    Let $q$ be an odd power of a prime $p$.
    The precise set of possible embedding degrees for simple supersingular
    abelian surfaces over $\fq$ is given in the following table.
    \begin{center}
    \begin{tabular}{|l|l|} \hline
    $p$ & Possible embedding degrees $k$ \\
    \hline
    2 & $\{ 1, 3, 6, 12 \} $ \\
    3 & $\{ 1, 3, 4 \} $ \\
    5 & $\{ 1, 3, 4, 5, 6 \} $ \\
    $\ge 7$ & $\{ 1, 3, 4, 6 \} $ \\
    \hline
    \end{tabular}
    \end{center}
\end{thm}

We note that other embedding degrees, such as $k=2$,
may be realised using non-simple abelian surfaces.
Since large embeddings degrees are of the most interest,
we focus on the cases where $k$ is $4$, $5$, $6$ and $12$.

\begin{description}
\item[$k=4$]
    The CM curve $y^2 = x^5 + A$ over $\F_p$
    where $p>2$ and $p \equiv 2, 3 \pmod{5}$
    is supersingular, and has embedding degree $4$.
    More generally,
    reductions of the CM curves listed by van Wamelen \cite{VanWamelen1}
    have embedding degree~$4$.
    These curves are discussed in Sections~\ref{Sec-3} and \ref{Sec-VW}.

\item[$k=5$]
    The curves $y^2 = x^5 - x \pm 1$
    where $p = 5$,
    described by Duursma and Sakurai \cite{DS},
    are supersingular and have embedding degree $5$.
    We discuss these curves in Section~\ref{sec:k5}.

\item[$k=6$]
    An abelian variety over $\fq$ has embedding degree $6$
    if its characteristic polynomial of Frobenius
    is of the form $T^4 - q T^2 + q^2$.
    A result of Howe, Maisner, Nart and Ritzenthaler \cite[Theorem 1]{HMNR}
    implies that such abelian varieties have a principal polarisation
    if and only if $p \not\equiv 1 \pmod{3}$.
    Hence,
    Jacobians of curves of genus~$2$
    can have embedding degree~$6$ only when $p \not\equiv 1\pmod{3}$.
    In Section~\ref{sec:twist},
    we give an algorithm to construct
    supersingular curves with embedding degree $6$
    when $p \equiv 2 \pmod{3}$ and $p \ge 5$,
    by taking suitable twists of the curve $y^2 = x^6 + 1$.

\item[$k=12$]
    The curves $y^2 + y = x^5 + x^3 + b$ over $\F_{2^m}$
    where $b=0,1$ are supersingular,
    with embedding degree~$12$.
    We consider these curves in Section~\ref{sec:char2}.
\end{description}

\section{Curves with embedding degree $4$: CM induced by an automorphism}
\label{Sec-3}

In this section,
we work with $q=p$ such that $p > 2$
and $p \equiv 2, 3 \pmod{5}$.
Consider the curve $\tilde{\C}$ defined over $\Q$ by
\[
    \tilde{\C}\colon y^2=x^5 + 1 .
\]
The curve $\tilde{\C}$ has an automorphism $\rho_5$ of order $5$
defined by
$$
    \rho_5: (x,y)\longmapsto (\zeta_5 x, y) ,
$$
where $\zeta_5$ is a primitive fifth root of unity over $\Q$.
The automorphism $\rho_5$
induces an endomorphism of $\Jac(\tilde{\C})$,
which we also denote~$\rho_5$.
The minimal polynomial of $\rho_5$ is the same as that of $\zeta_5$,
so $\End^0(\Jac(\tilde{\C}))$
contains the CM-field $\Q(\zeta_5)$.

Reducing $\tilde{\C}$ modulo $p$,
we obtain a curve $\C$ defined over $\fp$.
Since $p \not\equiv 1 \pmod 5$,
the endomorphism $\rho_5$ reduces
to a non-trivial endomorphism of $\Jac (\C)$,
also denoted~$\rho_5$.
This endomorphism was first used as a distortion map
by Choie and Lee~\cite{ChoieLee}.

\begin{rem}
    If $p \equiv 1 \pmod{5}$,
    then $\Jac(\C)$ is ordinary.
    If $p \equiv 4 \pmod{5}$,
    then $\Jac(\C)$ is supersingular but not simple.
    This explains our restriction to
    $p \equiv 2, 3 \pmod 5$.
\end{rem}

\begin{lem}
    The Jacobian $\Jac( \C )$ is
    $\fp$-simple,
    supersingular,
    and has embedding degree $4$.
\end{lem}

\begin{proof}
    Observe that $5$ does not divide $p-1$,
    so for each value of $y$ in $\F_{p}$
    there is a unique value $x = (y^2 - 1)^{1/5}$
    yielding a point $(x,y)$ in $\C(\F_p)$.
    Since $\C$ has a single point at infinity,
    we have $\#\C( \F_p ) = p+1$.
    Similarly,
    since $5$ does not divide $p^2  -1$,
    we obtain
    $\#\C( \F_{p^2} ) = p^2 +1$.
    It follows that
    $\#\Jac( \C )(\F_p ) = p^2 + 1$,
    and that
    the characteristic polynomial of the $p$-power
    Frobenius endomorphism $\pi$ on $\Jac (\C)$ is $P(T) = T^4 + p^2$.
    This polynomial is irreducible over $\Z$,
    so $\Jac( \C )$ is simple (but not absolutely simple).
    We may also deduce from the form of $P(T)$
    that $\C$ is supersingular (see \cite{S-X,Gal}).

    It remains to compute the embedding degree. 
    If $r$ is an odd prime 
    dividing $\#\Jac( \C )(\F_p )$,
    then $r$ divides $p^2 + 1$;
    hence
    $r$ divides $p^4-1$,
    and does not divide $p^i-1$ for any $i$ less than $4$.
    We conclude that $\Jac(\C)$ has embedding degree $4$.
\end{proof}

Our goal is to show that
for any pair of divisor classes on $\C$,
there is a suitable distortion map of the form $\pi^i \rho_5^j$
for some $i$ and $j$.
The first step towards establishing this result
is to show that the $\Q$-algebra $\End^0( \Jac(\C ))$
is generated as a $\Q$-module
by maps of the form $\pi^i \rho_5^j$.

Consider the non-commutative subring
$\Z[\rho_5,\pi]$ of $\End(\Jac(\C))$
generated by $\rho_5$ and $\pi$.
We let
$\Q[ \rho_5, \pi ]$ denote the non-commutative $\Q$-algebra
$\Z[ \rho_5, \pi ] \otimes_{\Z} \Q $.
Since
$\Z[ \rho_5, \pi ]$ is a finitely generated $\Z$-module,
$\Q[ \rho_5, \pi]$ is a finite dimensional $\Q$-vector space.
Note that
since the characteristic polynomial of $\pi$
has nonzero constant term,
there exists an element $\pi^{-1}$ of $\Q[\rho_5,\pi]$
such that $\pi^{-1}\pi = \pi\pi^{-1} = 1$.


\begin{lem}
\label{lemm-34}
    Let $\pi$ and $\rho_5$ be as above.
    Then $\pi^j \rho_5 \pi^{-j} = \rho_5^{(p^j)}$
    for all $j \ge 0$.
\end{lem}
\begin{proof}
    Clearly
    $
        \pi^j\rho_5(x,y)
        = (\zeta_5^{(p^j)} x^{(p^j)}, y^{(p^j)})
        = \rho_5^{(p^j)}\pi^j(x,y)
    $.
\end{proof}

Lemma~\ref{lemm-34}
implies that 
the inner automorphism
$\phi \mapsto \pi\phi\pi^{-1}$
of $\Q[ \rho_5, \pi ]$
has order~$4$ and fixes $\Q(\rho_5)$,
and so
corresponds to a generator $\sigma$ of the cyclic group
$\Gal (\Q(\rho_5)/\Q)$.
Let 
$\rho_5^{\sigma^j}$ denote the map $\sigma^j(\rho_5)$.
Since $\Gal (\Q(\rho_5)/\Q)$ is cyclic of order $4$,
the maps $\rho_5^{\sigma^j} $ are distinct
for $0\le j \le 3$.

\begin{prop}
\label{prop}
    As a $\Q$-vector space,
    $\Q[\rho_5,\pi]$
    has a direct sum decomposition
    $$
        \Q[\rho_5, \pi]
        =
        \Q(\rho_5) \oplus \pi\Q(\rho_5) \oplus \pi^2\Q(\rho_5)
            \oplus \pi^3\Q(\rho_5)
        .
    $$
\end{prop}

\begin{proof}
    We will prove by induction
    that
    the sum $\bigoplus_{i=0}^{t}\pi^{i}\Q(\rho_5)$
    is direct
    for each $0 \le t \le 3$.
    For $t=0$ there is nothing
    to prove.
    For the inductive step,
    assume
    $U_{n}=\bigoplus_{i=0}^{n}\pi^{i}\Q(\rho_5)$ is direct
    for $0 \le n \le 2$;
    we will show that $U_{n} \cap\pi^{n+1}\Q(\rho_5) = \{ 0 \}$.
    Suppose the contrary:
    then there is a non-zero $z$ in $\Q(\rho_5)$
    such that $\pi^{n+1}z$ is in $U_{n}$.
    Dividing by $z$,
    we can write $\pi^{n+1} =z_0+\pi z_1+ \dots +\pi^{n} z_{n}$,
    with coefficients $z_i$ in $\Q(\rho_5)$
    for $0 \le i \le n$,
    and with at least one of the $z_i$ not zero.
    Let $\sigma$ be a generator of $ \Gal( \Q( \rho_5) / \Q )$
    satisfying $\rho_5^{\sigma} = \rho_5^p$.
    Lemma \ref{lemm-34}
    implies $\rho_5^{\sigma^j} \pi = \pi \rho_5^{\sigma^{j-1}}$,
    and thus $\rho_5^{\sigma^{(n+1)}}\pi^{n+1} = \pi^{n+1}\rho_5$.
    Hence
    $$
    \begin{array}{r@{\;=\;}l}
        0
        &
        \rho_5^{\sigma^{n+1}} \pi^{n+1} - \pi^{n+1} \rho_5
        \\ &
        \rho_5^{\sigma^{n+1}} \left( z_0+  \pi z_1+ \dots + \pi^n z_n \right) - \left( z_0 \rho_5 + \pi z_1 \rho_5 + \dots + \pi^n z_n\rho_5 \right)
        \\ &
        z_0 \rho_5^{\sigma^{n+1}} +  \pi z_1 \rho_5^{\sigma^{n}} + \dots + \pi^n z_n\rho_5^{\sigma^{1}} - z_0 \rho_5 - \pi z_1 \rho_5 - \dots -\pi^n z_n\rho_5
        \\ &
        z_0 (\rho_5^{\sigma^{n+1}} - \rho_5) + \pi z_1 (\rho_5^{\sigma^{n}} - \rho_5 )+ \dots +\pi^n z_n(\rho_5^{\sigma} -  \rho_5 )
        .
    \end{array}
    $$
    But $U_n$ is a direct sum,
    and $\rho_5^{\sigma^j} \ne \rho_5$ for $1 \le j \le 3$;
    hence $z_0 = z_1 = \dots = z_t = 0$, which is a contradiction.
\end{proof}

\begin{cor}
\label{cor-1}
    We have
    $$
        \End^0(\Jac(\C))
        =
        \Q[\rho_5, \pi]
        =
        \left\{
            \sum_{0\le i, j \le 3} \lambda_{i,j} \pi^i \rho_5^j
                : \lambda_{i,j} \in \Q
        \right\}
        .
    $$
\end{cor}

\begin{proof}
    We know $\Q(\rho_5)$ is a $4$-dimensional $\Q$-vector space,
    so by Proposition~\ref{prop}
    $\Q[\rho_5,\pi]$
    is a $16$-dimensional $\Q$-vector subspace of $\End^0(\Jac(\C))$.
    But
    $\End^0(\Jac(\C))$
    is itself $16$-dimensional,
    so $\Q[\rho_5,\pi] = \End^0(\Jac(\C))$.
    The second equality then follows
    on noting that $\{ 1, \rho_5, \rho_5^2, \rho_5^3 \}$
    is a $\Q$-basis for $\Q(\rho_5)$.
\end{proof}

Theorem~\ref{thm-exist}
implies the existence of a distortion map $\phi$
for every pair $(D_1,D_2)$ of non-trivial points of order $r$ on $\Jac(C)$:
that is, an endomorphism $\phi$ such that $e_r( D_1, \phi(D_2)) \neq 1$.
Now $\End( \Jac( \C ))$
is an order in $\Q[ \rho_5, \pi ]$
containing $\Z[ \rho_5, \pi ]$,
so by
Corollary~\ref{cor-1}
there exist rational numbers $\lambda_{i,j}$
such that
$\phi = \sum_{i,j} \lambda_{i,j} \pi^i \rho_5^j$
in $\Q[\rho_5,\pi]$.
Let $m$ denote the least common multiple of the
denominators of the $\lambda_{i,j}$;
note that the endomorphism $m \phi$ is an element of $\Z[ \rho_5, \pi ]$.

\begin{assumption}
\label{assumption-1}
    We assume that $\phi$ may be chosen such that $\gcd( m, r ) = 1$,
    where $\phi$, $m$ and $r$ are defined as above.
\end{assumption}

\begin{rem}
Assumption~\ref{assumption-1} holds
if $\Z[ \rho_5, \pi]$ is ``most'' of $\End( \Jac( \C ))$,
and seems to hold in practical examples.
However,
we have not proven that it is always satisfied
for the curves under consideration.
It is instructive to consider Assumption~\ref{assumption-1}
in the case where $\Jac(\C)$ is a supersingular elliptic curve $E$.
In this case,
$\End(E)$ is a maximal order $\O$
in a quaternion algebra $B = \Q[ \pi, \psi ]$
where $\pi$ is the $q$-power Frobenius
and $\psi$ is some other endomorphism.
Note that
$\alpha \O \alpha^{-1}$ is a maximal order in $B$
for every $\alpha$ in $B$;
so maximal orders in $B$ can be very far
from $\Z[ \pi, \psi]$.
Thus Assumption~\ref{assumption-1}
may not be true in general.
However,
following the arguments of \cite{GR},
we may suppose that $E$ has been constructed by the CM method,
in which case $\psi^2 = -d$ for some relatively small positive integer $d$.
Hence $\Z[ \pi, \psi ]$ is contained in $\O$,
and Assumption~\ref{assumption-1}
holds when $r$ is larger than $d$.
\end{rem}

\begin{thm}
\label{easy-case-thm}
    If Assumption~\ref{assumption-1} holds,
    then for
    all pairs $(D_1, D_2)$ of non-zero divisor classes on $\C$ of order $r$
    and
    all non-degenerate pairings $e_r$
    there exists a distortion map of the form $\pi^i \rho_5^j$
    with $0 \le i, j \le 3$.
\end{thm}

\begin{proof}
    Theorem~\ref{thm-exist} shows that there exists an endomorphism $\phi$
    that is a suitable distortion map for $(D_1,D_2)$;
    Corollary~\ref{cor-1} shows that $\phi$ is in $\Q[ \rho_5, \pi ]$.
    Under Assumption~\ref{assumption-1},
    we may take an integer $m$ prime to $r$ such that
    $m\phi$ is in $\Z[ \rho_5, \pi ]$
    and
    \[
        e_r( D_1,  m\phi( D_2 ) )
        = e_r( D_1, \phi( D_2))^m \neq 1
        ;
    \]
    so $m \phi$ is also a distortion map for $(D_1,D_2)$.
    Since $m\phi$ is an integer combination of the $\pi^i \rho_5^j$,
    we must have $e_r( D_1, \pi^i \rho_5^j( D_2 )) \neq 1$
    for some $0 \le i,j \le 3$
    (otherwise,
    if all $e_r( D_1, \pi^i \rho_5^j( D_2 )) = 1$,
    then $e_r( D_1,  m\phi( D_2 ) ) = 1$ by the linearity of the pairing).
\end{proof}

\begin{rem}
    Alternatively, one could use
    maps of the form $\rho_5^i \pi^j$
    in Theorem~\ref{easy-case-thm}.
\end{rem}

\begin{example}
    Let $D_1$ be a nonzero element of $\Jac(\C)[r]$
    defined over $\F_p$;
    note that $\pi( D_1 ) = D_1$.
    It is easy to show that, under Assumption~\ref{assumption-1},
    $e_r( D_1, \rho_5^j( D_1 )) \neq 1$
    for some $1 \le j \le 3$.
    This supports the suggestion in \cite{ChoieLee}
    of using $\rho_5$ as a distortion map.
    When implementing pairings,
    it is desirable to utilise denominator elimination
    to improve efficiency;
    to this end, the map $\rho_5^j$
    might be combined with a trace operation
    (see Scott \cite{Scott} for an example of
    this in the elliptic case).
\end{example}

\begin{rem}
    The results in this section easily generalise to the
    twists $y^2 = x^5 + A$ of $C$ (for nonzero $A$),
    and even more generally to the curves
    $y^2 = x^{2n+1} + A$ over $\fp$,
    where $2n+1$ is prime and $p$ is a
    primitive root modulo $2n+1$.
\end{rem}

\section{Curves with embedding degree $4$: Other CM curves}
\label{Sec-VW}

In \cite{VanWamelen1} and \cite{VanWamelen2},
van Wamelen
describes the $19$ isomorphism classes of curves of genus $2$ over $\Q$
whose Jacobians have CM by the ring of integers of a CM-field.
For each isomorphism class,
van Wamelen provides a representative curve $\tilde\C_i$ defined over $\Q$,
the CM-field $F_i := \End^0(\Jac(\tilde\C_i)$,
and an explicit partial description of
an endomorphism $\alpha_i$ of $\Jac(\tilde\C_i)$
such that $F_i = \Q(\alpha_i)$,
giving the ($x$-coordinates of) the image under $\alpha_i$
of the image of a generic point $(x,y)$ of $\tilde\C_i$ in $\Jac(\tilde\C_i)$.
One can derive a full description of the endomorphism $\alpha_i$
from this information (see Pujol\`as~\cite{pujolas} for details).

The curve $\tilde C$ of Section~\ref{Sec-3}
is a representative of the isomorphism class
corresponding to the CM-field $\Q(\zeta_5)$
in van Wamelen's tables.
In this section,
we generalize our treatment of $\tilde C$
to the other CM curves $\tilde C_i$.
Reducing each $\tilde{\C_i}$ modulo suitable inert primes $p$,
we obtain curves $\C_i$ over $\fp$
whose Jacobians are simple,
supersingular,
and whose characteristic polynomial of Frobenius is equal to $T^4 + p^2$.
These Jacobians therefore have a very similar
endomorphism structure to that of the Jacobian in Section~\ref{Sec-3}.

%

%
%

Let $\tilde\alpha_i$ be an endomorphism of $\Jac( \tilde{\C_i} )$
such that
$\End^0(\Jac( \tilde{\C_i} )) \cong \Q( \tilde\alpha_i )$;
the endomorphism supplied by van Wamelen suffices.
Note that $\tilde\alpha_i$ is defined over the quartic field $F_i$,
which has cyclic Galois group over $\Q$
for all of the curves in~\cite{VanWamelen1}.
Let $\alpha_i$ denote the image of $\tilde\alpha_i$ in $\End( \Jac( \C_i ))$.
Since $p$ is inert in the cyclic quartic field $F_i$,
it follows that $\alpha_i$  is defined over $\F_{p^4}$.

If $\pi$ is the $p$-power
Frobenius map then, as before, we have
\[
    \pi \alpha_i \pi^{-1} = \alpha_i^{(p)}
\]
where $\alpha_i^{(p)}$ denotes the map obtained from $\alpha_i$
by applying the $p$-power Frobenius
to the coefficients of $\alpha_i$.
It follows that 
the inner automorphism $\phi \mapsto \pi\phi\pi^{-1}$
generates $\Gal( \Q( \alpha_i )/\Q )$.
We may therefore prove an analogue of Proposition~\ref{prop}
for each curve $\C_i$.

\begin{prop}
    The non-commutative $\Q$-algebra $\Q[ \alpha_i, \pi ]$
    generated by $\alpha_i$ and $\pi$
    is a $16$-dimensional $\Q$-vector space,
    and (as $\Q$-vector spaces)
    $$
        \Q[ \alpha_i, \pi ]
        =
        \Q(\alpha_i) \oplus \pi\Q(\alpha_i) \oplus \pi^2\Q(\alpha_i)
            \oplus \pi^3\Q(\alpha_i)
        .
    $$
\end{prop}


As a result,
under the appropriate analogue of Assumption~\ref{assumption-1},
we may choose a distortion map of the form
$\pi^u \alpha_i^v$ with $0 \le u, v \le 3$
for any pair of
elements of $\Jac(C_i)[r]$.
The van Wamelen curves are therefore
suitable for cryptography,
in the sense that
one can easily find a distortion map
for every pair of divisors.

\begin{rem}
    In practice,
    evaluating the maps $\alpha_i$ of~\cite{VanWamelen2}
    is relatively complicated,
    making the distortion maps of the curves in this section
    relatively inefficient
    compared with those of
    the CM curve $y^2 = x^5 + 1$ described in Section~\ref{Sec-3}.
\end{rem}

\begin{rem}
    One could also construct curves with distortion maps
    by reducing CM curves defined over number fields other than $\Q$.
\end{rem}

\section{Curves with embedding degree $5$}
\label{sec:k5}

The curves $C : y^2 = x^p - x + b$ over $\F_p$
with $b = \pm 1$
have been studied by Duursma and Sakurai \cite{DS},
and efficient pairing computation on these curves
was studied by Duursma and Lee \cite{DL}.
Our interest is in the genus $2$ case,
so in this section
we consider the curves
$$
    C: y^2 = x^5 - x + b
$$
over $\F_q$,
where $q=5^m$ for some $m$ coprime to $10$,
and $b = \pm 1$.
The distortion map proposed by Duursma and Lee \cite{DL} is
$$
        \psi(x,y) = ( \rho - x, 2 y ),
$$
where $\rho$ is an element of $\F_{5^5}$ such that $\rho^5 - \rho + 2b = 0$.

The characteristic polynomial of the ($q$-power) Frobenius for these curves is
\[
    P^{\pm}_m(T)
    =
    T^4 \pm 5^{(m+1)/2} T^3 + 3 \cdot  5^m T^2 \pm 5^{(3m+1)/2} T + 5^{2m}
    .
\]
Observe that
\[
        P^{+}_m(T) P^{-}_m(T)
    =
    T^8 + q T^6 + q^2 T^4 + q^3 T^2 + q^4
    ,
\]
and hence that $(T^2 - q) P^{+}_m(T) P^{-}_m(T) = T^{10} - q^5$.
Let $N := \#\Jac( C )(\F_{q} )$.
Since
$N$ is equal to either $P^+_m(1)$ or $P^-_m(1)$,
it follows that
$N$ divides $q^5 -1$;
hence the embedding degree is $k=5$
for large prime-order subgroups of $\Jac(C)(\F_q)$.
Note that the characteristic polynomial of
the $q^5$-power Frobenius is $(T^2 - 5^{5m})^2$,
and that
the full $N$-torsion is defined over $\F_{q^{10}}$ but not over $\F_{q^5}$.
Since $\Jac(C)(\F_{q^k})[N] \cong (\Z / N\Z)^2$,
this case is as easily handled as
the elliptic curve case in Section~\ref{sec:intro}.

\begin{lem}
\label{lem61}
    If $D_1$ and $D_2$ are non-zero points of prime order $r$
    in $\Jac(C)(\F_q)$,
    then $\psi$ is a suitable distortion map
    with respect to the Tate pairing.
\end{lem}
\begin{proof}
    As in the elliptic curve case,
    the pairing of $D_1$ with $D_2$
    is defined over $\F_q$,
    and is therefore trivial.
    On the other hand,
    $\psi(D_2)$ is a non-zero $r$-torsion divisor
    which is not defined over $\F_q$.
    Thus $\{ D_1, \psi(D_2)\}$
    is a basis for $\Jac(C)(\F_{q^{k}})[r]$,
    and by non-degeneracy of the Tate pairing
    we have $e_r( D_1, \psi( D_2 )) \ne 1$.
\end{proof}

\begin{rem}
    Note that Lemma~\ref{lem61}
    is only stated for the Tate pairing.
    This is because while the Tate pairing is known to be
    non-degenerate for points defined over the field
    $\F_q( \mu_r ) = \F_{q^5}$,
    we are only guaranteed that the Weil pairing is non-degenerate
    when working over
    $\F_q(\Jac(C)[ r ]) = \F_{q^{10}}$.
\end{rem}

If $C$ is a curve defined over $\F_q$,
then for any integer $n$
there is a natural homomorphism
$\Tr : \Jac(C)(\F_q^n) \to \Jac(C)(\F_q)$
defined by
$$
    \Tr(D) := \sum_{i=0}^{n-1}\pi^i(D) .
$$
This map is called the {\em trace} (of Frobenius);
its kernel is called the {\em trace-zero subgroup} of $\Jac(C)(\F_q)$.

\begin{lem}
\label{lem62}
    If $D$ is a point in $\Jac(C)(\F_q)$,
    then $\psi(D)$ lies in the trace-zero subgroup
    of $\Jac(C)(\F_{q^5})$.
\end{lem}

\begin{proof}
    We first prove the result for divisors of the form
    $D = (P) - (\infty)$
    where $P$ is a point in $C( \F_q )$.
    Observe that the $y$-coordinate of $\psi( P )$
    is in $\F_q$.
    Hence the $y$-coordinates of $\pi^j( \psi( P ))$ for $0\le j\le 4$
    are all equal,
    say to some $y_0$,
    while the $x$-coordinates
    are all distinct.
    The function
    $(y - y_0)$ on $C$
    therefore has divisor equal to
    $(\psi(P)) + (\pi( \psi(P))) + \cdots + (\pi^4(\psi(P))) - 5(\infty)$,
    which is the trace of $(\psi(P)) - (\infty) = \psi(D)$.
    Hence, the trace of $\psi(D)$ is zero.

    The same argument applies to divisors of the form
    $(P_1) + (P_2) - 2(\infty)$ where $P_1$ and $P_2$ are in $C( \F_q )$.
    A similar argument applies when
    $P_1$ and $P_2$ are Galois conjugates in $C( \F_{q^2} )$:
    the $y$-coordinates of the $\pi^j( \psi( P_i))$
    take either the same value $10$ times,
    or two different values $5$ times each.
\end{proof}


Using Lemmas~\ref{lem61} and \ref{lem62},
it immediately follows (as in \cite{Verheul2,GR})
that $\psi$ and combinations of $\pi^j$ 
are sufficient as distortion maps
for all pairs of points of order $r$ in $\Jac( C)(\F_{q^5} )$.
On the other hand,
since $\psi^2=-1$,
it is clear that
$\Q[ \pi, \psi ]$ is a $\Q$-algebra of dimension $8$;
so combinations of $\pi$ and $\psi$
are not sufficient to act as distortion maps
for all pairs of points in $\Jac(C)[r]$.

To obtain generators for the full endomorphism ring we use the fact
that $C$ is isomorphic over $\F_{5^5}$ to the curve $C' : Y^2 = X^5
- X$ (the isomorphism is $\eta (x,y) = (x - \alpha, y)$ where
$\alpha \in \F_5$ satisfies $\alpha^5 - \alpha + b = 0$).  On $C'$
we have the automorphism $\phi'(X,Y) = (2X, \beta Y)$ where $\beta
\in \F_{5^2}$ satisfies $\beta^2 = 2$. Defining $\phi = \eta^{-1}
\phi' \eta$ gives the automorphism
\[
   \phi( x,y) = (2x - \alpha, \beta y)
\]
on $C$.  Since $\phi$ is not defined over $\F_{5^{5}}$ it follows
that $\phi$ does not lie in $\Q[ \pi, \psi ]$ and hence $\Q[ \pi,
\psi, \phi ] = \End^0( \Jac(C))$.  It follows that if the Weil
pairing is used then a distortion map of the form $\pi^u \psi^v
\phi^w$ with $0 \le u \le 3, 0 \le v, w \le 1$ may be used.

\begin{rem}
    Since the embedding degree is odd,
    the fact that the image of the distortion map is a trace zero divisor
    does not mean that
    the usual denominator elimination techniques
    for pairing implementation may be used.
    However,
    the ate pairing approach
    (see \cite{HSV,GHOTV})
    may be used to obtain a very simple pairing algorithm,
    with no final exponentiation required.
    This does not imply that characteristic $5$ curves
    are fast for pairing-based cryptography,
    since $5$ is not a very natural base
    for computer implementation of finite field arithmetic.
\end{rem}

\section{Curves with embedding degree $6$}
\label{sec:twist}

Let $p$ be an odd prime such that $p \equiv 2 \pmod{3}$.
Let $\zeta_6$ be
a primitive sixth root of unity over $\F_p$,
and set $\zeta_3:=\zeta_6^2$.

We wish to construct a curve $C'/\F_p$
with embedding degree~$6$.
The characteristic polynomial of Frobenius on $\Jac(C')$
must therefore be
$T^4 - pT^2 + p^2$.
Following \cite[p.32]{HNR},
we obtain $C'$
by twisting the curve $C : y^2 = x^6 + 1$
with respect to its automorphism
$u : (x,y) \mapsto (\zeta_3/x,y/x^3)$.
To find a defining equation for $C'$,
we need to find an isomorphism $\phi : C' \to C$
defined over $\overline{\F}_p$
such that $\phi^{(p)} \circ \phi^{-1} = u$.
We can assume $\phi$ is of the form
$$
        \phi(X,Y)
    = (x,y)
    = \left( \frac{aX + b}{cX+d} , \frac{Y}{(cX + d)^3} \right)
    ,
$$
with $ad - bc \not= 0$;
the curve $C'$ will then have a defining equation
$$
    C' : Y^2 = (a X + b)^6 + (c X + d)^6.
$$
We need $\phi^{(p)} = u\circ\phi$
--- that is,
$$
    \left(\frac{a^pX + b^p}{c^pX + d^p},\frac{Y}{(c^pX + d^p)^3}\right)
    =
    \left(\frac{\zeta_3cX + \zeta_3d}{aX + b},\frac{Y}{(aX + b)^3}\right)
    .
$$
To find particular solutions for $a$, $b$, $c$ and $d$
we begin by setting $a = c^p$ and $b = d^p$.
Now $a^p = \zeta_3 c$ and $b^p = \zeta_3 d$,
so we need values for $c$ and $d$ such that
$$
    c^{p^2 - 1} = d^{p^2 - 1} = \zeta_3 .
$$
Therefore, we choose some $\gamma$ in $\overline{\F}_{p}$ satisfying
$\gamma^{p^2 - 1} = \zeta_3$, and set $c = \gamma$. One can show
that $\gamma \in \F_{p^6}$. Note that $(c/d)^{p^2 - 1} = 1$, so
$c/d$ is an element of $\F_{p^2}$. On the other hand, we know $ad
\not= bc$, so $c^pd \not= d^pc$: hence $c/d$ is not an element of
$\F_p$. Therefore, setting $d = \zeta_3 c$, we obtain a solution
$$
    a = \gamma^p,\quad
    b = \zeta_3^{-1} \gamma^p,\quad
    c = \gamma, \text{ and }
    d = \zeta_3\gamma
    .
$$

We want to find suitable distortion maps for $\Jac(C')$.
We could proceed as in Section~\ref{Sec-3},
and obtain a conditional result
depending on the final assumption about the denominators.
However,
since $C'$ is isomorphic to $C$,
we can consider the problem of
finding distortion maps for $\Jac(C)$ instead:
indeed,
if $A \subset \End(\Jac(C))$
is a suitable set of distortion maps for $\Jac(C)$,
then $\phi^{-1} A \phi$ will be
a suitable set of distortion maps for $\Jac(C')$.
This approach allows us to take advantage
of the splitting behaviour of $\Jac(C)$,
and thus to obtain an unconditional result.


Let $E$ be the elliptic curve defined over $\F_p$
by $E : y^2=x^3+1$,
and let $\pi_E$ denote the $p$-power Frobenius on $E$.
The curve $E$ has an automorphism $\rho_3$
defined over $\F_{p^2}$ by
$(x,y) \mapsto (\zeta_3 x, y)$.
\begin{lem}
    With $E$, $\rho_3$ and $\pi_E$ as above,
    \begin{enumerate}
    \item   $E$ is supersingular,
    \item   the characteristic polynomial of $\pi_E$ is $X^2 + p$, and
    \item   $\Z[\pi_E,\rho_3]$ is an order of index $3$ in $\End(E)$.
    \end{enumerate}
\end{lem}
\begin{proof}
    One easily checks that $\pi_E\circ\rho_3 \not= \rho_3\circ\pi_E$
    when $p \equiv 2 \pmod{3}$,
    so $\End(E)$ is non-commutative;
    it follows that $E$ is supersingular.
    The characteristic polynomial of $\pi_E$
    has the form $T^2 - tT + p$,
    where $-2\sqrt{p} \le t \le 2\sqrt{p}$;
    but since $E$ is supersingular,
    $p$ divides $t$~\cite[Theorem~V.3.1]{Silverman},
    and the only such $t$ is $0$.
    Hence the characteristic polynomial of Frobenius is $T^2 + p$.
    Since $E$ is supersingular,
    $\End(E)$ is isomorphic to a maximal order of the quaternion algebra
    ramified at $p$ and $\infty$;
    its discriminant is therefore $p$~\cite[Corollary~5.3]{Vigneras}.
    Explicit calculation shows that
    $\Z[\pi_E,\rho_3]$ is an order of discriminant $3p$,
    and thus an order of index $3$ in $\End(E)$.
\end{proof}

Let $f : C \to E$
(resp.~$f' : C \to E$)
be the morphism defined by $f(x,y)= (x^2,y)$
(resp.~$f'(x,y)=(1/x^2,y/x^3)$).
We define homomorphisms
$$
\begin{array}{r@{\;}r@{\;}c@{\;}l}
    	\mu : 	& E \times E 	& \longrightarrow & 
         \Jac(C) \\
		& (P, Q)	& \longmapsto	  & f^*(P) + {f'}^*(Q)
                - 4 P_{\infty}\\
\end{array}
$$
and
$$
\begin{array}{r@{\;}r@{\;}c@{\;}l}
    	\tilde\mu: & \Jac(C) & \longrightarrow & E \times E \\
		& P+Q-2 P_{\infty} & \longmapsto & (f_*(P)+f_*(Q),f'_*(P)+f'_*(Q))  \\
\end{array}$$
where $P_{\infty}=(0,1,0) \in C$.

Observe that
$\tilde{\mu} \circ \mu = [2]_{E \times E}$
and
$\mu \circ \tilde{\mu} = [2]_{\Jac(C)}$,
so $\mu$ and $\tilde\mu$ are $(2,2)$-isogenies.
We can therefore define an injective (group) homomorphism
$$
\begin{array}{r@{\;}c@{\;}l}
    T :
    \End(\Jac(C))
    &
    \longrightarrow
    &
    \End(E \times E)
    \\ 
    \psi
    &
    \longmapsto
    &
    \tilde{\mu} \circ \psi \circ \mu
    .
    \\
\end{array}
$$
While $T$ is not a ring homomorphism,
one easily checks that
$$T(\psi) T(\psi')= 2 T(\psi \psi')$$
for all endomorphisms $\psi$ and $\psi'$ of $\Jac(C)$.

Let $\chi$ and $\rho_6$ be
the automorphisms of $C$ defined by
$$
    \chi(x,y) = ( 1/x ,  y/x^3 )
    \text{\quad and\quad }
    \rho_6(x,y) = ( \zeta_6 x, y ) ;
$$
we use the same notations for the induced endomorphisms of
$\Jac(C)$. Note that $\pi_C,\chi,\rho_6$ all preserve $P_{\infty}$.
Let $\pi_C$ denote the $p$-power Frobenius on $C$,
and let $A=\Z[\pi_C,\chi,\rho_6]$ be
the subring of $\End(\Jac(C))$
generated by $\pi_C$, $\chi$, and $\rho_6$.
We will compute an upper bound
for the index of $T(A)$ in $\End(E \times E)=M_2(\End(E))$.
It suffices to compute
the images of $\pi_C$, $\chi$, and $\rho_6$. 

\begin{lem}
\label{endo-images}
    The images of $\pi_C$, $\chi$, and $\rho_6$ in $\End(E\times E)$
    are given by
    $$
        T(\pi_C) =
        2 \left(\begin{array}{cc}
                \pi_E & 0 \\
                0 & \pi_E \\
            \end{array}
        \right),
        \
        T(\chi) =
        2 \left(\begin{array}{cc}
            0  & 1 \\
            1 & 0\\
        \end{array}\right),
        \ \text{and}\
        T(\rho_6) =
        2 \left(\begin{array}{cc}
            \rho_3 & 0 \\
            0  & - \rho_3^2 \\
        \end{array} \right)
        .
    $$
\end{lem}
\begin{proof}
    We show this for $\rho_6$, the other computations being
    similar. We will abuse notations by dropping the $P_{\infty}$s
    since they remain unchanged in the computations.
    Consider a point $(x^2,y)$ on $E$.
    We have $f^*(x^2,y)=(x,y) + (-x,y)$ in $\Jac(C)$,
    so
    $$
        \rho_6 \circ f^*(x^2,y) = (\zeta_6 x,y)+ (-\zeta_6 x, y).
    $$
    Now
    $$
    \begin{array}{r@{\;=\;}l}
        f_* \circ \rho_6 \circ f^*(x^2,y)
        &
        ((\zeta_6 x)^2,y)+((-\zeta_6 x)^2,y) \\
        &
        2 (\zeta_3 x^2,y) \\
        &
        ([2] \circ \rho_3) (x^2,y)
    \end{array}
    $$
    whereas
    $$
    \begin{array}{r@{\;=\;}l}
        (f'_* \circ \rho_6 \circ f^*)(x^2,y)
        &
        \left( 1/(\zeta_6 x)^2, y/(\zeta_6 x)^3 \right)
        +
        \left(1/(-\zeta_6 x)^2, y/(-\zeta_6 x)^3 \right)
        \\
        &
        \left(1/(\zeta_3 x),-y/x^3\right)
        +
        \left (1/(\zeta_3 x),y/x^3\right)
        \\
        & 0
        .
    \end{array}
    $$
    In the same way ,
    $
        {f'}^*\left(1/x^2,y/x^3\right)
        =
        (x,y) + (-x,-y)
    $,
    so
    $$
    \begin{array}{r@{\;=\;}l}
        (f'_* \circ \rho_6 \circ {f'}^*)\left(1/x^2,y/x^3\right)
        &
        \left(1/(\zeta_6 x)^2,y/(\zeta^6 x)^3\right)
        +
        \left(1/(-\zeta_6 x)^2,-y/(-\zeta^6 x)^3\right)
        \\
        &
        \left(1/(\zeta_3 x^2),-y/x^3\right)
        +
        \left(1/(\zeta_3 x^2),-y/x^3\right)
        \\
        &
        ([-2] \circ \rho_3^2) \left(1/x^2,y/x^3\right)
    \end{array}
    $$
    whereas
    $
        (f_* \circ \rho_6 \circ {f'}^*)\left(1/x^2,y/x^3\right)
        = 0
    $.
\end{proof}

Using Lemma~\ref{endo-images},
we see that the images of
$1+\rho_6^3$, $1-\rho_6^3$, $\chi+\chi \rho_6^3$ and $\chi-\chi \rho_6^3$
in $\End(E\times E)$
are the ``projectors''
$$
    4 \left(\begin{array}{cc}
        1 & 0 \\
        0 & 0 \\
    \end{array}\right),
    \
    4 \left(\begin{array}{cc}
        0  & 0 \\
        0 & 1\\
    \end{array}\right),
    \
    4 \left(\begin{array}{cc}
        0 & 1 \\
        0  & 0 \\
    \end{array} \right),
    \
    \text{ and }
    \
    4 \left(\begin{array}{cc}
        0 & 0 \\
        1  & 0 \\
    \end{array} \right).
$$
Composing with the images of $\pi_C$ and $\rho_6$,
we get
$$
    M_2(4 \Z[\pi_E,\rho_3])
    \subset
    T(A)
    \subset(M_2(\End(E)).
$$
Since $\Z[\pi_E,\rho_3]$ has index $3$ in $\End(E)$,
the index of $T(A)$ in $M_2(\End(E))$
divides $2^8\cdot3^4$.
Now,
suppose $r$ is a prime different from $p$, $2$ and $3$.
Following the proof of Theorem~\ref{thm-exist},
by tensoring with $\Z/r\Z$ we get an isomorphism
$$
    T_r :
    \End(\Jac(C)) \otimes \Z/r\Z
    \stackrel{\sim}{\longrightarrow}
    M_2(\End(E)) \otimes \Z/r\Z \simeq M_4(\Z/r\Z)
$$
and $T_r(A)$ is of index dividing $2^8\cdot3^4$ in $M_4(\Z/r\Z)$.
Thus,
if $D_1$ and $D_2$ are non-zero elements of $\Jac(C)[r]$,
we can find a map $\Phi$ in $M_4(\Z/r\Z)$
such that $e_r(D_1,\Phi(D_2)) \neq 1$.
Then $2^8\cdot3^4 \Phi = T_r(\psi)$
for some $\psi$ in $A$,
and
$$
    e_r(D_1,\psi(D_2))
    =
    e_r(D_1,[2^8\cdot3^4]\Phi(D_2))
    =
    e_r(D_1,\Phi(D_2))^{2^8\cdot3^4}
    \neq
    1 .
$$
We have proven the following theorem.

\begin{thm}
    Let $r$ be a prime different from $2$, $3$ and $p$.
    For all pairs $D_1$, $D_2$ of non-zero elements of $\Jac(C')[r]$,
    there exists a suitable distortion map in the ring
    $\phi^{-1} \Z[\pi_C,\chi,\rho_6] \phi$.
\end{thm}

\begin{rem}
    Using the construction of \cite[p.~12]{HLP}
    with the group-scheme isomorphism $\eta : E[2] \to E[2]$
    mapping $(-1,0)$ to itself and $(\zeta_6,0)$ to $(1/\zeta_6,0)$,
    we have $\Jac(C) \simeq (E \times E)/\textrm{Graph}(\eta)$.
    Moreover,
    if $\lambda$ is the canonical polarization on $\Jac(C)$
    and $\lambda_{E \times E}$ the split polarization on $E \times E$,
    then $\lambda_{E \times E} = \hat{\mu} (2 \lambda) \mu$.
    Thus if $D_1$ and $D_2$ are elements of $\Jac(C)[r]$,
    then
    \begin{eqnarray*}
        e^{\lambda}_r(D_1,D_2)^8
        &=&
        e^{2 \lambda}_r(2 D_1,2 D_2)
        =
        e^{2 \lambda}_r(\mu \tilde{\mu} D_1,\mu \tilde{\mu} D_2)
        \\
        &=&
        e^{\lambda_{E \times E}}_r(\tilde{\mu} D_1,\tilde{\mu} D_2)
        \\
        &=&
        e^{\lambda_E}_r(f_*(D_1),f_*(D_2))
            \cdot e^{\lambda_E}_r(f'_*(D_1),f'_*(D_2))
        .
    \end{eqnarray*}
    In particular,
    pulling back two divisors on $\Jac(C')$
    first to $\Jac(C)$ and then to $E \times E$,
    we see that the computation of the pairing on $\Jac(C')$
    is in fact equivalent to
    the computation of twice the pairing on the elliptic curve $E$.
\end{rem}

\section{Curves with embedding degree~$12$}
\label{sec:char2}

The curves $\C: y^2 + y = x^5 + x^3 + b$ over $\F_{2^m}$ with
$b=0,1$ were studied by van der Geer and van der Vlugt in
\cite{VanderGeer1} and \cite{VanderGeer2}, in view of their
applications to coding theory. Throughout this section, the ground
field is $\F_{2^m}$, where $m \equiv \pm 1 \pmod 6$.

The characteristic polynomial of
the Frobenius endomorphism $\pi$ over $\F_{2^m}$ is
$$
    P_m^{\pm}(T)
    =
    T^4\pm2^{(m+1)/2}T^3+2^mT^2\pm2^{(3m+1)/2}T+2^{2m}
    ,
$$
so $\Jac(\C)$ is supersingular and simple over $\F_{2^m}$.
We have
$$
    P_m^+(T)P_m^-(T)
    =
    T^8-2^{2m}T^4+2^{4m}
$$
and
$$
    (T^8-2^{4m})(T^8+2^{2m}T^4+2^{4m})P_m^+(T)P_m^-(T)
    =
    T^{24}-2^{12m}
    ,
$$
so the embedding degree is $k=12$.

Theorem~\ref{thm-exist}
shows that a distortion map $\phi$ exists
for every pair $(D_1,D_2)$
of points in $\Jac(\C)[r]$.
We will now give a set of maps that contains a distortion map for
any pair of divisors in $\Jac(\C)(\F_{2^{12m}})[r]$,
by exhibiting a basis of $\End^0(\Jac (\C))$.

The automorphisms of $\C$ are of the form
$$
    \sigma_{\omega}:\,(x,y)\longmapsto(x+{\omega},y+s_2x^2+s_1x+s_0)
$$
where ${\omega}$ is any root of the polynomial
$$
\begin{array}{r@{\;}l}
    &
    x^{16}+x^8+x^2+x \\
    = &
    (x^6+x^5+x^3+x^2+1)(x^3+x^2+1)(x^3+x+1)(x^2+x+1)(x+1)x
    ,
\end{array}
$$
and where $s_2={\omega}^8+{\omega}^4+{\omega}$,
$s_1={\omega}^4+{\omega}^2$,
and $s_0$ is a root of $y^2 + y ={\omega}^5+{\omega}^3$
(note that $s_0+1$ is the other root).
For each $\omega$,
we arbitrarily fix one of the corresponding $s_0$,
and denote the resulting automorphism $\sigma_\omega$.
We interpret these automorphisms as elements of $\End(\Jac(\C))$.
One can verify that they satisfy the relations
$$
    \sigma_{\omega}\sigma_{\omega'}
    =
    \pm\sigma_{\omega'}\sigma_{\omega}
    =
    \pm \sigma_{\omega+\omega'}
    .
$$

Fix a root $\tau$ in $F_{2^{6}}$ of $x^6+x^5+x^3+x^2+1$,
and set
$\xi=\tau^4+\tau^2$,
$\rho=\tau^2+\tau+1$, and
$\theta=\tau^4+\tau^2+\tau$.
Note that $\xi$, $\rho$ and $\theta$ are
roots of the cubic and quadratic factors above.
We have
$\theta^2=\theta+1$,
$\tau^8=\tau+1$, and
$\theta+\tau=\xi$.
As before,
we let $\Z[ \pi, \sigma_\tau, \sigma_\theta ]$
denote the non-commutative ring generated by $\pi$ and $\sigma_\tau$,
and write
$
    \Q[\pi,\sigma_{\tau}, \sigma_\theta ]
    =
    \Z[\pi,\sigma_{\tau}, \sigma_\theta ]\otimes \Q
$
for the algebra generated by $\pi$, $\sigma_\tau$ and $\sigma_\theta$.

\begin{prop}
    The $\Q$-algebra $\Q[\pi,\sigma_{\tau},\sigma_\theta]$
    is a $16$-dimensional $\Q$-vector space
    with a direct sum decomposition
    $$
        \Q[\pi,\sigma_{\tau}, \sigma_\theta ]
        =
        \Q(\pi) \oplus \sigma_{\tau}\Q(\pi) \oplus \sigma_{\theta}\Q(\pi)
            \oplus \sigma_{\xi}\Q(\pi)
        .
    $$
    Furthermore,
    $\End^0(\Jac(\C))=\Q[\pi,\sigma_{\tau},\sigma_\theta ]$.
\end{prop}

\begin{proof}
    Let $F = \Q( \pi )$;
    note that $F$ is a $4$-dimensional $\Q$-vector space.
    One easily checks that the relations
    $$
    \begin{array}{l}
        \pi \sigma_{\omega} = \pm \sigma_{{\omega}^{2^m}} \pi,
        \\
        \sigma_{\tau}^2 = -1,
        \\
        \pi^3 \sigma_{\tau} \pi^{-3}
            = \pm \sigma_{\tau^{2^3}}
            = \pm \sigma_{\tau + 1}
            = \pm \sigma_{\tau} \sigma_1,
        \text{ and }\\
        \pi \sigma_\theta \pi^{-1}
            = \pm \sigma_{\theta^2}
            = \pm \sigma_{\theta + 1}
            = \pm \sigma_1 \sigma_\theta
    \end{array}
    $$
    hold,
    where $\sigma_1$ is the automorphism $( x,y) \mapsto (x+1, y + x^2)$.

    Let $A := F \oplus\sigma_{\tau} F $;
    the sum is direct,
    because $\sigma_{\tau}$ is not in $F$.
    We see that $A$ is an $8$-dimensional $\Q$-vector space.
    Note that $A$ is not an algebra.

    We now show that $A \oplus \sigma_{\xi} F $ is direct.
    Assume the contrary:
    that is, that there is some non-zero $z$ in $F$
    such that $\sigma_{\xi} z$ lies in $A$.
    Dividing by $z$,
    we must have
    \[
            \sigma_{\xi} = z_1 + \sigma_{\tau} z_2
    \]
    for some $z_1$ and $z_2$ in $F$.
    Now, since $\xi \in \F_{2^3}$
    we have $\sigma_{\xi} = \pi^3 \sigma_{\xi} \pi^{-3}$.
    Using the relations above,
    we see
    $$
    \begin{array}{r@{\;=\;}l}
        z_1 + \sigma_{\tau} z_2
        &
        \pi^3 \sigma_{\xi} \pi^{-3}
        \\ &
        \pi^3 (z_1 + \sigma_{\tau} z_2 ) \pi^{-3}
        \\ & 
	z_1 \pm \sigma_{\tau^{2^3}} z_2
        \\ & 
	z_1 \pm \sigma_{\tau} \sigma_1 z_2
        .
    \end{array}
    $$
    Since $A$ is a direct sum and $\sigma_1 \ne \pm 1$,
    we have $z_2 = 0$:
    that is,
    $\sigma_{\xi}$ must lie in $F$,
    which is a contradiction
    since $\sigma_{\xi}$ does not commute with $\pi$.

    Finally,
    we show that $(A \oplus \sigma_{\xi} F ) \oplus \sigma_{\theta} F$
    is direct.
    Assuming the contrary,
    we have
    \[
            \sigma_{\theta} = z_1 + \sigma_{\tau} z_2 +  \sigma_{\xi} z_3
    \]
    for some $z_1$, $z_2$ and $z_3$ in $F$.
    Again using the relations above, we have
    $$
    \begin{array}{r@{\;=\;}l}
        0 &
        \sigma_1 \sigma_{\theta} \pi^3 \pm \pi^3 \sigma_{\theta}
        \\ &
	\sigma_1( z_1 + \sigma_{\tau} z_2 + \sigma_{\xi} z_3 )\pi^3
                \pm \pi^3 (z_1 + \sigma_{\tau} z_2 + \sigma_{\xi} z_3 )
        \\ & 
	\sigma_1( z_1 + \sigma_{\tau} z_2 + \sigma_{\xi} z_3 )\pi^3
                \pm (z_1 \pm \sigma_{\tau^{2^3}} z_2 \pm \sigma_{\xi^{2^3}} z_3 ) \pi^3
        \\ & 
	(\sigma_1 \pm 1) z_1\pi^3
        \pm (\sigma_1 \pm \sigma_{1} )\sigma_{\tau} z_2 \pi^3
        \pm (\sigma_1 \pm 1) \sigma_{\xi} z_3 \pi^3.
    \end{array}
    $$
    Since $A \oplus \sigma_{\xi} F$ is direct
    and $\sigma \not= \pm 1$,
    we must have $z_1 = z_3 = 0$.
    Therefore,
    $\sigma_{\theta} = \sigma_{\tau} z$ for some $z$ in $F$;
    but this is a contradiction,
    since $\sigma_{\theta}$ is defined over $\F_{2^2}$,
    while $\sigma_{\tau}$ is defined over $\F_{2^6}$.
    We conclude that
    $(A \oplus \sigma_{\xi} F ) \oplus \sigma_{\theta} F$
    is direct.

    Thus
    $
        \Q[\pi,\sigma_{\tau},\sigma_{\theta}]
        =
        F \oplus \sigma_{\tau}F \oplus \sigma_{\xi}F
            \oplus \sigma_{\theta}F
    $,
    and is therefore a $16$-dimensional $\Q$-vector space.
    Since $\End^0(\Jac(\C))$ is $16$-dimensional
    and contains $\Q(\pi,\sigma_{\tau},\sigma_{\theta})$,
    we have
    $ \End^0(\Jac(\C)) = \Q[\pi,\sigma_{\tau},\sigma_{\theta}] $.
\end{proof}

Our claim that a distortion map for any pair of divisors
may be chosen from
the maps $\pi$, $\sigma_{\tau}$, and $\sigma_{\theta}$
follows.
Indeed,
for any pair of points $(D_1,D_2)$ in $\Jac(\C)[r]$,
we may choose a distortion map $\phi$
as a $\Q$-linear combination of the endomorphisms
$\pi^i$, $\pi^j\sigma_\tau$, $\pi^k\sigma_\theta$ and $\pi^l\sigma\xi$.
If we assume that we may choose coefficients
such that the least common multiple $m$
of their denominators
is coprime to $r$,
then
$m\phi$ lies in $\Z[\pi,\sigma_\tau,\sigma_\theta]$,
and is a suitable distortion map.

\section{Conclusions and future work}

We have given several examples of distortion maps
for supersingular Jacobians of genus $2$ curves
with embedding degree $4$, $5$, $6$ and $12$.
We have proven,
subject to a reasonable assumption,
that these maps are sufficient for all applications.

One natural problem for future study
is to show that the assumption holds
for the curves considered in this paper.
Another problem is
to consider similar problems in the genus $3$ case,
although Rubin and Silverberg \cite{RS} have shown that
there is little motivation for using high genus curves in pairing applications.

\section*{Acknowledgements}

The authors would like to thank Ryuichi Harasawa, David Kohel, and Enric~Nart.

The work described in this paper has been supported in part by the
European Commission through the IST Programme under Contract IST-2002-507932
ECRYPT.  The information in this document reflects only the author's views, is
provided as is and no guarantee or warranty is given that the information is
fit for any particular purpose.  The user thereof uses the information at its
sole risk and liability.

This research was also supported by the
EU SOCRATES/ERASMUS programme
and by the EPSRC.

\end{document}